\documentclass[11pt, reqno]{amsart}
\usepackage{amsmath, amssymb, amsthm, graphicx,marvosym}
\usepackage{cancel}

\usepackage{tikz}
\usepackage{amsmath}
\usepackage{amssymb}
\usepackage{amsthm}
\usepackage{enumerate}
\newtheorem{thm}{Theorem}
\newtheorem{lemma}{Lemma}[section]

\newtheorem{definition}[lemma]{Definition}

\makeatletter \@addtoreset{equation}{section} \makeatother

\newcommand{\la}{\lambda}
\newcommand{\RR}{\mathbb{R}}
\newcommand{\beq}{\begin{equation}}
\newcommand{\eeq}{\end{equation}}

\def\com#1{\quad{\textrm{#1}}\quad}

\def\nn{\nonumber}

\def\com#1{\quad{\textrm{#1}}\quad}

\def\nn{\nonumber}

\def\({\left(\begin{array}{cccccc}}
\def\){\end{array}\right)}

\def\bes{\begin{eqnarray}}
\def\ees{\end{eqnarray}}

\begin{document}

\title[Existence and regularity in wave equations]
{Existence and regularity of solutions in nonlinear wave equations}

\author{Geng Chen}
\address{Geng Chen, School of Mathematics,
Georgia Institute of Technology, GA, USA
30332. 
({\tt email: gchen73@math.gatech.edu}).}
\author{Yannan Shen}
\address{Yannan Shen, Department of Mathematical Sciences, 
University of Texas at Dallas, Richardson, TX, USA 75080.
({\tt email:  yxs135630@utdallas.edu})}

\begin{abstract}
In this paper, we study the global existence and regularity of 
H\"older continuous solutions for a series of nonlinear partial differential equations describing nonlinear waves.
\end{abstract}

\maketitle
\tableofcontents



\section{Introduction}

In this paper, we consider the existence and regularity of weak solutions for two families of nonlinear wave equations parameterized by $\lambda$
\beq\label{spe}
	u_{tx}+f'(u)\, u_{xx}+\lambda\, f''(u)\, (u_x)^2=0\,, 
\eeq
and 
\beq\label{vwl}
u_{tt}-c^2(u)\, u_{xx}-2\lambda\, c(u)\, c'(u)\, (u_x)^2=0\,,
\eeq
with constant parameter
\[ 
0\leq\lambda\leq1\,.
\]
Here $x\in\mathbb{R}$ is the spatial variable and $t\in\mathbb{R}^+$ is the time variable. The wave speed $c(u)>0$. 
Such equations can be formally written as
\beq\label{spe2}
u_{tx}+(f'(u))^{1-\lambda}\, \bigl[(f'(u))^{\lambda}\,u_{x}\bigr]_x=0\,, 
\eeq
\beq\label{vwl2}
u_{tt}-c^{2-2\la}(u)\, [c^{2\la}(u)u_{x}]_x=0\,,
\eeq
which include several important and interesting models when $\lambda$ takes different values. 

%
\begin{itemize}
\item For equation (\ref{spe}) with\vspace{.2cm}

\begin{itemize}
\item[$\la=1$:] Scalar hyperbolic conservation law.\vspace{.2cm}
\item[$\la=\frac{1}{2}$:] An equation considered in \cite{BZZ}. When $f(u)=\frac{1}{2}u^2$,  equation (\ref{spe}) is Hunter-Saxton equation modeling nematic liquid crystal \cite{BC,BHR, HR,HS,HZ95a,HZ95b}.\vspace{.2cm}
\item[$\la=0$:] A wave equation in unitary direction, $u_{tx}+f'(u)u_{xx}=0$. \end{itemize}\bigskip


\item For equation (\ref{vwl}) with\vspace{.2cm}

\begin{itemize}

\item[$\la=1$:] Wave equation modeling elasticity
\[u_{tt}-\bigl(F(u)\bigr)_{xx}=0,\quad \text{with} \quad F(u)=\int c^2(u)du\,,\] 
or 
isentropic Euler equations in Lagrangian coordinates, also called p-system:
	\begin{eqnarray}
	u_t-\omega_x=0\label{ps}\\
	\omega_t-F_x=0\nn
	\end{eqnarray}
	with $\omega=\int u_t~dx$\,, and $F$ denotes pressure. See \cite{Dafermos} for details.
	\vspace{.2cm}
\item[$\la=\frac{1}{2}$:] Variational wave equation modeling nematic liquid crystal \cite{BZ,HR,ghz,ZZ05a}. \beq\label{vw}u_{tt}-c(u)(c(u)u_x)_x=0\,.\eeq
\vspace{.1cm}
\item[$\la=0$:]  A nonlinear wave equation \[u_{tt}-c^2(u)u_{xx}=0\] which is the one dimensional case of \beq\label{vwlin}u_{tt}-c^2(u)\Delta u=0\,,\eeq which was studied in \cite{lindblad}.
	
\end{itemize}
\end{itemize}

\bigskip
%

One common feature of these systems is the finite time gradient blowup of solutions even with smooth initial data, when $0<\la\leq 1$, \cite{CPZ,ghz,HS,lax0}. The motivation why we connect all these equations together is to understand the variation of regularity for weak solutions of these wave equations as $\la$ changes. The first appearance of Equation (\ref{vwl}) was in \cite{ghz2}.
%

We notice that the regularities of weak solutions for equations in the form of (\ref{spe}) (or (\ref{vwl})) with $\la=0$, $\frac{1}{2}$ and $1$, respectively, are totally different. We summarize the existing results on regularity of weak solutions for these three cases in the following table.
\begin{itemize}
\item When $\la=1$, (\ref{spe}) and (\ref{vwl})
can be written in the form of hyperbolic conservation laws. It is well known that solutions in these equations in general have discontinuities (shock waves) even when initial data are smooth, c.f. \cite{G3,CPZ,G5,G8,Dafermos,lax0}. BV existence for solution of (\ref{ps}) with small amplitude is available in \cite{LG}.
\vspace{.1cm}

\item When $\la=\frac{1}{2}$, solutions for Hunter-Saxton equation and  variational wave equation 
are in general only H\"older continuous with exponent $1/2$ because of the possible gradient blowup, c.f. \cite{BZ,BZZ,HZ95a} for global existence and  \cite{GCZ,ghz,HS} for gradient blowup. \vspace{.1cm}

\item When $\la=0$, there is still no global existence for classical solutions available for general large initial data. But 
we tend to expect that the solution in this case has better regularity than solutions in previous two cases, because of the global-in-time existence for classical radially symmetric small solution in \cite{lindblad} and  the study in this paper for large data solution.
\end{itemize}

By the discussion for three cases with $\la=0$, $\frac{1}{2}$ and $1$, it is very tentative for us to guess that the solution for  (\ref{spe}) or (\ref{vwl}) is more regular when $\la$ is decreasing. More intuitively,
as $\la$ decreases, i.e. ``more" $f'(u)$ or $c(u)$ comes out of the bracket in \eqref{spe2} or \eqref{vwl2}, we conject that weak solution of \eqref{spe} or \eqref{vwl} has better 
regularity. 

In this paper, we partially prove this conjecture. 
First, we show that the conjecture is true for (\ref{spe}) with $\la\in(0,\frac{1}{2}]$ by constructing weak solutions whose
regularities vary on $\la$. Especially, when  $\la\in(0,\frac{1}{3}]\cup\frac{1}{2}$, the solution is H\"older continuous 
on both $x$ and $t$ with exponent $1-\la$. 

Secondly, we provide some numeric evidences supporting that the solution for wave equation (\ref{vwl}) with $\la\in(0,\frac{1}{3}]$ is H\"older continuous with exponent
$1-\la$ when gradient blowup happens. Especially, very loosely speaking, when $\la$ is very close to zero,  
H\"older space with exponent $1-\la$ is getting ``close" to the space of $C^1$ functions.

In fact, for (\ref{vwl}), we construct
a semi-linear system, then the problem whether solution $u(x,t)$ is H\"older continuous with exponent
$1-\la$ is changed to an equivalent problem whether variables $p$ and $q$ defined in (\ref{0pq_def}) in the semi-linear system are bounded away from zero and infinity, which can be more easily tested by numerical methods than the first problem. We do numerical experiments on several examples, in all of which $p$ and $q$ are indeed bounded away from zero and infinity, although gradient blowup happens in finite time.

We expect this work can help unveiling the mystery in H\"older continuous solutions for quasi-linear hyperbolic systems. To our limit of knowledge, studies on H\"older continuous solutions in wave equations are still very limited. Especially the classification of equations whose solutions are H\"older continuous with different exponents is wide open.

In next two subsections, we introduce the main results for (\ref{spe}) and (\ref{vwl}).

%
\subsection{{Global existence for  (\ref{spe}) with $\la\in(0,\frac{1}{2}]$}}
For (\ref{spe}), we focus on the case when $\la\in(0,\frac{1}{2}]$. 
The energy law for (\ref{spe}) for the smooth solution
is
\beq\label{spe_energy}
\Big(|u_x|^\frac{1}{\la}\Big)_t+\Big(f'(u)|u_x|^\frac{1}{\la}\Big)_x=0\,.
\eeq


In this paper, 
we consider the initial boundary value problems for (\ref{spe}) on the region
$(x,t)\in \mathbb{R}^+\times\mathbb{R}^+$ when $\la\in(0,\frac{1}{2}]$, with initial and boundary conditions
		\beq\label{ID1} 
			u(0,t)=0,\qquad u(x,0)=:{u}_{0}(x)\in W_{loc}^{1,\frac{1}{\la}}( \RR^+)\,,
		\eeq
		and  a compatibility condition
		\beq\label{CC}
			u_0(0)=0 \com{and} u'_0(0)=0\,.
		\eeq
Here $W_{loc}^{1,\frac{1}{\la}}( \RR^+)$ is the Sobolev space with standard notation.
				
Throughout this paper, we assume that $f(u)$ is a $C^2$ function satisfying
\beq\label{f_assp}
f'(0)\geq 0,\qquad |f''(u_1)-f''(u_2)|\leq L|u_1-u_2|\,,\qquad \forall\ u_1,\,u_2\in \mathbb{R}
\eeq
for some constant $L$, and $|f''(u)|$ is uniformly bounded above. 
The assumption that $f'(0)\geq 0$ protects that the wave on the boundary $x=0$ does not flow in an outward direction.

We first define the weak solution.
\begin{definition}{\bf\text (Weak solution)}\label{def}
The function $u(x,t)$, defined for all $(x,t)\in\mathbb{R}^+\times\mathbb{R}^+$,
is a weak solution for \eqref{spe}, if initial and boundary conditions \eqref{ID1} and \eqref{CC} 
are satisfied pointwisely and
	\begin{itemize}
	\item[i.] The equation (\ref{spe}) is satisfied in the weak sense		
	\beq\label{weak}
			\int_{0}^\infty\int_{0}^\infty \Big\{-u_{x}\, 
			(\phi_t+f'(u)\,\phi_x)+(\lambda-1)\, f''(u)u_x^2\,\phi\Big\}\, dx\, dt=0,
		\eeq
		for any test function $\phi\in C_c^1(\mathbb{R}^+\times\mathbb{R}^+)$.
	\item[ii.] For any fixed $t>0$, the function $u(\cdot,t)$ is in $W_{loc}^{1,\frac{1}{\la}}( \RR^+)$, hence is locally H\"older continuous with exponent $1-\la$ by the Sobolev embedding Theorem. 
	\end{itemize}			
\end{definition}
\smallskip

Then we give the main theorem in this paper for (\ref{spe}).

{\begin{thm}\label{main}
The initial boundary value problem
\eqref{spe}\eqref{ID1}\eqref{CC} with $\la\in(0, \frac{1}{2}]$ exists a weak solution $u(x,t)$ under Definition \ref{def}.  When $\la\in(0,\frac{1}{3}]\cup\frac{1}{2}$,
$u(x,t)$ is locally H\"older continuous on both $x$ and $t$ with exponent $1-\la$. 
\end{thm}}

To prove Theorem \ref{main}, inspired by the energy dependent characteristic coordinates introduced in \cite{BC2} for Camassa-Holm equation and in \cite{BZ}  for variational wave equation, we introduce an independent variable $Y$, which  dilates the possible gradient blowup due to the concentration of characteristics.
Then we establish a semi-linear system for some unknowns on independent variables $Y$ and time. One crucial unknown is 
\[
\xi=\frac{(1+u_x^2)^\frac{1}{2\la}}{Y_x}\,,
\]
where $Y_x$ measures the dilation rate of characteristics.
By showing $\xi$ is bounded, which can be roughly understood as that the dilation rate of characteristics $Y_x$ is balanced by the energy density, we can prove the global existence 
for the semi-linear system when $\la\in(0,\frac{1}{2}]$. Finally, using an inverse transformation,
we construct the weak solution for (\ref{spe}).
When $\la\in(0,\frac{1}{3}]\cup\frac{1}{2}$, the solution has better regularity.

%
\subsection{Analysis on nonlinear wave equation \eqref{vwl}}

When $\lambda=\frac{1}{2}$, for \eqref{vwl}, 
energy conservative H\"older continuous solutions with exponent $1/2$ has been proved in \cite{BZ} by introducing new characteristic variables as independent variables. This method has also
been used for a series of variational wave equations in nematic liquid crystal \cite{CZZ12,GCZ,ZZ10,ZZ11}.
In these results and also our paper, the wave speed $c(u)$ is assumed to be uniformly positive and bounded.

In this paper, inspired by (\ref{spe}), we derive a semi-linear system for \eqref{vwl}. 
Although we still cannot conclude any global existence results, we expect
this system could create a framework in proving the global existence of H\"older continuous solutions for \eqref{vwl} in the future. 

When $\la\in(0,\frac{1}{3}]$, the only issue left towards the global existence of H\"older continuous solution with exponent $1-\la$ is that we cannot find the uniform $L^\infty$ bound on two variables $p$ and $q$ defined in \eqref{0pq_def}, which take similar role as  $\xi$ for (\ref{spe}). 
Hence, 
the issue whether the solution $u$ is 
H\"older continuous with exponent $1-\la$ is changed to another issue whether variables $p$ and $q$ in the semi-linear system are bounded away from zero and infinity, at the breakdown of classical solution.
Especially, if we use numeric method to study these two issues, the latter one is much simpler than the first one.

Several numeric experiments are given in this paper, which all indicate that $p$ and $q$ are bounded away from zero and infinity even when
classical solution of \eqref{vwl} breaks down. As a consequence, it is reasonable to expect that the solution $u(x,t)$ is H\"older continuous with exponent $1-\la$ even when the gradient blowup happens in finite time, where $\la\in(0,\frac{1}{3}]$.
\bigskip

The rest of the paper is divided into two sections. In Section 2, we consider system \eqref{spe} and prove Theorem \ref{main}. In Section 3, we will discuss the wave equation \eqref{vwl}.

\section{Wave in a unitary direction}
In this section, we consider system \eqref{spe} and prove Theorem \ref{main}. 

In Subsection \ref{subsec_new}, we first define a new coordinate $(Y,T)$. Then based on \eqref{spe}  and initial condition \eqref{ID1}\eqref{CC}, we 
derive a semi-linear system for several unknowns on new independent variables $Y$ and $T$ . 

In Subsection \eqref{subsec_exist}, we prove existence and uniqueness of solution for the new semi-linear system on $(Y,T)$ coordinates. 
 
Finally, in  Subsection \eqref{subsec_reverse}, after making an inverse transformation on the constructed solution on $(Y,T)$ coordinates, 
we recover a weak solution for (\ref{spe}) on the $(x,t)$-coordinates and complete the proof of Theorem \ref{main}.
\subsection{New coordinates\label{subsec_new}}
In this subsection, we derive some equations valid for smooth solutions of \eqref{spe}. We denote the variables
\beq\label{S_R_def}
\left\{
\begin{array}{l}
 S := u_t + f'(u) u_x\,\\
 R := u_x\,.
\end{array}\right.
\eeq
By (\ref{spe}), we have
\beq\label{RxSx}
\left\{
\begin{array}{rcl}
 S_x &=& (1-\lambda)\, f''(u)\, R^2\,\\
 R_t +f'(u) R_x &=& -\lambda f''(u)\, R^2\,.
\end{array}\right.
\eeq
In order to reduce the equation (\ref{spe}) into a semi-linear system, it is convenient to change the independent variables.
The equation of the characteristic is
\beq\label{char}
\frac{d x^c(t)}{d t}=f'\bigl(u(x^c(t),t)\bigr) .
\eeq
We denote the characteristic passing through the point $(x,t)$ as 
\[
a\mapsto\, x^c(a;\, x,\, t) \com{or equivalently} b\mapsto\, t^c(b;\, x,\, t),
\]
where $a$ and $b$ are the time and space variables of the characteristic, respectively.
\begin{figure}[htb]
\centering
\includegraphics[scale=.28]{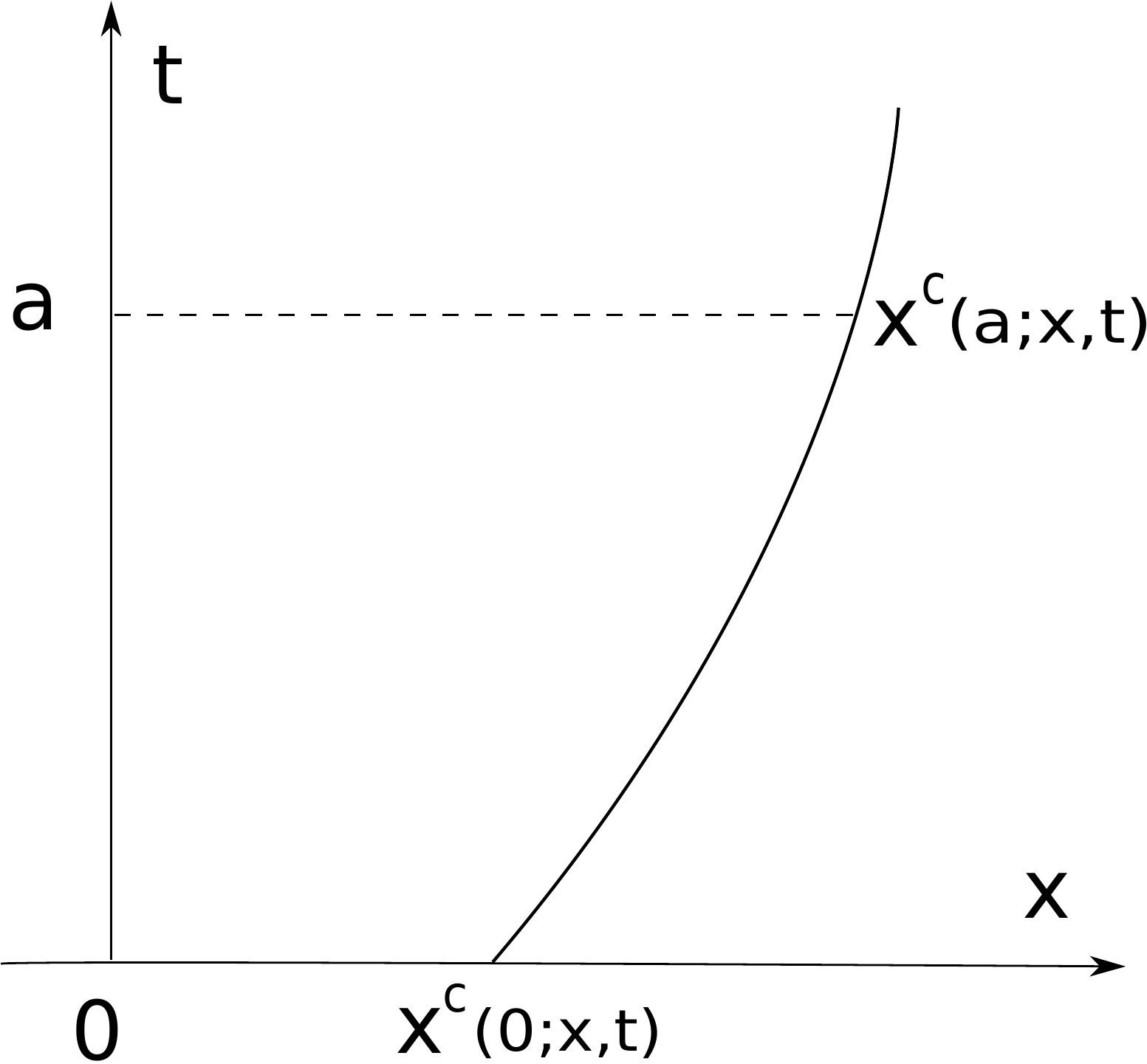}
\caption{A characteristic $a\mapsto\, x^c(a;\, x,\, t)$ passing $(x,t)$.}
\label{p1}
\end{figure}
Then we introduce new coordinates $(Y,T)$, such that
\beq\label{Y_def}
Y\equiv Y(x,t):=\left\{\begin{array}{l}
					\int_0^{x^c(0;\, x,\, t)} (1+ R^2(x',0))^\frac{1}{2\la}\,d x', \\
					\ \com{when the characteristic passing $(x,t)$ interacts $t=0$;}
					\\
					\\
					-t^c(0;\, x,\, t)\, f'(0)\\
					\ \com{when  the characteristic passing $(x,t)$ interacts $x=0$,}
				  \end{array}\right.
\eeq 
with $(x,t)\in{\mathbb R}^+\times{\mathbb R}^+$ and
\beq
T\equiv T(x,t):=t\,.
\eeq
Clearly, $Y$ is constant along a characteristic $x^c$ by its definition. So
\beq\label{Y_T_eqn}
Y_t +f'(u) Y_x=0, \quad T_t=1 \com{and} T_x= 0.
\eeq
Using (\ref{Y_T_eqn}), for any smooth function $m$, we have
\beq\left\{\begin{array}{l}\label{trans}
m_t +f'(u) m_x=m_Y\,(Y_t +f'(u) Y_x) +m_T\, (T_t +f'(u) T_x)=m_T,\\
m_x=m_Y\, Y_x+m_T\, T_x=m_Y\, Y_x.
\end{array}\right.
\eeq

Then we derive a semi-linear system on $(Y,T)$-coordinates.
In order to complete the system, we introduce several new variables: 
\beq\label{def_lgq}
v:=2\arctan{u_x} \com{and} \xi:=\frac{(1+R^2)^
\frac{1}{2\lambda}}{Y_x}.
\eeq
Hence
\beq\label{tri}
\frac{1}{1+R^2}=\cos^2 \frac{v}{2} \com{and}
 \frac{R}{1+R^2}=\frac{1}{2}\sin v.
\eeq
By (\ref{RxSx}) and (\ref{trans}), we have
\beq\label{u_eqn2}
u_Y=\frac{u_x}{Y_x}=\frac{1}{2}\xi\,\sin{v}\, (\cos^2{\frac{v}{2}})^{\frac{1}{2\la}-1},
\eeq
\beq\label{g_eqn}
v_T=\frac{2\,R_T}{1+R^2}=-2\la\, f''(u)\,\frac{R^2}{1+R^2}=-2\la\, f''(u)\sin^2{\frac{v}{2}},
\eeq
and also using (\ref{Y_T_eqn}), we have
\beq\label{q_eqn}
\left.
\begin{array}{rcl}
\xi_T&=&(\frac{(1+R^2)^\frac{1}{2\la}}{Y_x})_T\smallskip\\
	&=&\frac{1}{\la}\frac{(1+R^2)^{\frac{1}{2\la}-1}\,R\, R_T}{Y_x} -\frac{(1+R^2)^\frac{1}{2\la}}{Y_x^2}(Y_x)_T\smallskip\\
	&=&-\frac{f''(u)(1+R^2)^{\frac{1}{2\la}-1}R^3}{Y_x}- \frac{(1+R^2)^\frac{1}{2\la}}{Y_x^2}\bigl((Y_t+f'(u)\,Y_x)_x-f''(u)RY_x \bigr)\smallskip\\
	&=&f''(u)\frac{(1+R^2)^\frac{1}{2\la}}{Y_x}\,\frac{R}{1+R^2} \bigskip\\
	&=&\frac{1}{2}f''(u)\, \xi\, \sin v.
	\end{array}
	\right.
\eeq
Summarizing (\ref{u_eqn2})$\sim$(\ref{q_eqn}), we have
a semi-linear system:
\beq\left\{\begin{array}{rcl}\label{semi}
u_Y&=&\frac{1}{2}\xi\,\sin{v}\, (\cos^2{\frac{v}{2}})^{\frac{1}{2\la}-1}\,,\vspace{.2cm}\\
v_T&=&-2\la\, f''(u)\sin^2{\frac{v}{2}}\,,\vspace{.2cm}\\
\xi_T&=&\frac{1}{2}f''(u)\, \xi\, \sin v\,.
\end{array}\right.
\eeq

Furthermore, 
\beq\label{u_eqn}
u_T=u_t +f'(u)\, u_x=S,
\eeq
and
\beq\label{S_eqn}
S_Y= \frac{S_x}{Y_x}= (1-\lambda)\, f''(u)\,\xi\, \sin^2 \frac{v}{2}\,(\cos^2{\frac{v}{2}})^{\frac{1}{2\la}-1}\,.
\eeq
Hence we have another semi-linear system
\beq\left\{\begin{array}{rcl}\label{semi2}
u_T&=&S\,,\vspace{.2cm}\\
S_Y&=& (1-\lambda)\, f''(u)\,\xi\, \sin^2 \frac{v}{2}\,(\cos^2 {\frac{v}{2}})^{\frac{1}{2\la}-1}\,,\vspace{.2cm}\\
v_T&=&-2\la\, f''(u)\sin^2{\frac{v}{2}}\,,\vspace{.2cm}\\
\xi_T&=&\frac{1}{2}f''(u)\, \xi\, \sin v\,.
\end{array}\right.
\eeq

It is easy to see that semi-linear systems \eqref{semi} and \eqref{semi2} are both invariant under translation by $2\pi$ in $v$. 
It would be more precise to use $e^{iv}$ as variable. For simplicity, we use $v\in[-\pi,\pi]$ with endpoints identified.
\subsection{Existence on the new coordinates\label{subsec_exist}}
In this subsection, we prove the existence of solution for
(\ref{semi}) with initial and boundary data converted from
(\ref{ID1}) and (\ref{CC}). To avoid the confusion,
the reader should be aware that in this section we solve variables $(u,v,\xi)$ 
by system (\ref{semi}) when $\la\in(0,\frac{1}{2}]$, and we do not use the equations (\ref{spe}) and definition (\ref{def_lgq}) except when we assign the initial and boundary data. We also use  (\ref{semi2}) when $\la\in(0,\frac{1}{3}]\cup\frac{1}{2}$ to show better
regularity for the solution. We will recover a weak solution for (\ref{spe}) on the $(x,t)$-coordinates in subsection \eqref{subsec_reverse},

\subsubsection{The initial boundary value problem on new coordinates\label{3.1}}

The initial lines $T=0$ and $t=0$ are the same line. 
By (\ref{Y_def}), the curve $Y:=\Gamma_b(T)$ on the $(Y,T)$-plane transformed from $\{x=0,\ t \geq 0\}$ is
\beq\label{Y_def0}
Y=\Gamma_b(T)=-f'(0) T\,.
\eeq
Recall $f'(0)\geq0$.

After transformation from $(x,t)$-coordinates to $(Y,T)$-coordinates, 
set $(\RR^+,\RR^+)$ changes to a new set named $\Omega$:
\beq\label{O_def}
\Omega:=\{(Y,T);\ Y\geq\Gamma_b(T),\ T\geq 0\}.
\eeq



We consider initial boundary value problem of (\ref{semi}) on $\Omega$
with following initial and boundary data given by \eqref{ID1}.

The initial data on $(Y,0)$ with $Y\geq 0$ are
\beq \label{bc}
\left\{\begin{array}{rcl}
	u(Y,0)&:=&u_{0}\bigl(x(Y,0)\bigr),\\ 
	v(Y,0)&:=&2\arctan (u_{0}'\bigl(x(Y,0)\bigr)),\\ 
	\xi(Y,0)&:=&1.
\end{array}\right.
\eeq

The boundary conditions on $Y=\Gamma_b(T) $ with $T\geq 0$ are
\beq \label{bc1}
\left\{\begin{array}{rcl}
	u(\Gamma_b(T),T)&:=&0,\\ 
	v(\Gamma_b(T),T)&:=&0,\\ 
	\xi(\Gamma_b(T),T)&:=&1.
\end{array}\right.
\eeq 

%
%


{\begin{thm}\label{main3}
Assume all conditions on initial and boundary data in Theorem \ref{main} hold. Then the corresponding problem
(\ref{semi}) with boundary data (\ref{bc})(\ref{bc1}) has a unique solution defined
for all $(Y,T)\in\Omega$.

Moreover, one has the stability of the solution.
Assume that a sequence of $C^1$ functions
$u_{0}^k$  satisfy
\[ u_{0}^k\to u_{0}\,,\quad
 (u_{0}^k)_x\to (u_{0})_x\,,\quad
\]
uniformly on any bounded subset of ${\mathbb R}^+$. Then one has the
convergence of the corresponding solutions for (\ref{semi}):
$$(u^k,  v^k, \xi^k)\to (u, v, \xi)$$
uniformly on bounded subsets of $\Omega$.

\end{thm}}
\begin{proof}
The proof is based on the locally Lipschitz continuity of the right hand side of equations (\ref{semi}). Actually, on any bounded domain 
\[
	\Omega_r:=\left\{ (Y,T);\ Y\geq \Gamma_b(T),\ 0\leq T\leq r,\ Y\leq r \right\}.
\] 
Since $|f''(u)|$ is uniformly bounded above, using the last equation in (\ref{semi}), we could find a priori upper bound of $\xi$, i.e. $\xi<e^{\frac{1}{2}r\max_{u\in \Omega_r}(f''(u))}\,.
$
Then as  long as $\lambda\in(0, \frac{1}{2}]$, the right hand side of equations (\ref{semi}) is Lipschitz continuous on $(u,v,\xi)$ in $\Omega_r$.

Introduce a space of functions 
\[
	\Theta_r := \left\{ f : \Omega_r \mapsto \mathbb{R};\  \| f \|_*:= 
	{\text{ess}}\sup_{(Y,T)\in \Omega_r} e^{-\kappa(T+|Y|)} |f(Y,T)|<\infty \right\},
\]
where $\kappa$ is a suitably large constant. It is straightforward to construct a solution $(u,v,\xi)(Y,T)$ with $(Y,T)\in\Omega_r$ as a fixed point in $\Theta_r\times\Theta_r\times \Theta_r$, using the integral forms of (\ref{semi})
and the fact that the right hand side of (\ref{semi}) is Lipschitz. 
When $r\rightarrow\infty$, we can get a unique solution defined on the whole domain of $\Omega$. As a consequence of the fixed point argument,  solutions have stability in the sense of Theorem \ref{main3}. We refer  readers to \cite{BZ} for more details. 

\end{proof}

By Theorem \ref{main3}, we have the regularity of $(u,v,\xi)$.
In fact, the equation  (\ref{semi}) implies, on any bounded set of  $\Omega$,
when $\gamma\in(0,\frac{1}{2}]$:
\begin{itemize}
\item $u$ is Lipschitz continuous w.r.t $Y$, measurable  
w.r.t $T$.  
\item $v, \xi$ are Lipschitz continuous w.r.t $T$, measurable  
w.r.t $Y$.  
\item $u, v, \xi$ have finite $L^\infty$ norm, and $\xi>0$.
\end{itemize}
\bigskip

When $\gamma\in(0,\frac{1}{3}]\cup\frac{1}{2}$, right hand side of \eqref{semi2} is also locally Lipschitz continuous, hence we could also
prove global existence of the solutions using \eqref{semi2} by same argument. This solution is exactly the same solution of \eqref{semi} because
\[
u_{YT}=(1-\lambda)\, f''(u)\,\xi\, \sin^2 \frac{v}{2}\,(\cos^2 {\frac{v}{2}})^{\frac{1}{2\la}-1}=u_{TY}\,.
\]
Hence we have
\begin{itemize}
\item $u, S$ are Lipschitz continuous w.r.t $Y$, measurable  
w.r.t $T$.  
\item $u, v, \xi$ are Lipschitz continuous w.r.t $T$, measurable  
w.r.t $Y$.  
\item $u, v, S, \xi$ have finite $L^\infty$ norm, and $\xi>0$.
\end{itemize}

\subsection{Existence on $(x,t)$-coordinates\label{subsec_reverse}}


The map $(Y,T)\mapsto (x,t)$ from $\Omega$ to $(\RR^+,\,\RR^+)$ can be constructed by the following procedure.

First we use $t=T$ and
\beq\label{trans1}
x_Y=\xi(\cos^2\frac{v}{2})^{\frac{1}{2\la}} \com{and} 
x_T=f'(u)
\eeq
to do the inverse transformation from  $(Y,T)$ to $(x,t)$. 
It is easy to check that the two equations in (\ref{trans1}) are equivalent:
\[
x_{TY}=\frac{1}{2}f''(u)\xi\,\sin{v}\, (\cos^2{\frac{v}{2}})^{\frac{1}{2\la}-1}=x_{YT}\,.
\]
So we can recover the function $x(Y,T)$ by integrating either
$x_Y=\xi(\cos^2\frac{v}{2})^{\frac{1}{2\la}}$ or
$x_T=f'(u)$\,. 
By (\ref{trans1}), it is easy to recover  (\ref{Y_T_eqn})  and (\ref{trans}).

For any smooth function $m$, 
\beq\label{trans2}
\xi(\cos^2\frac{v}{2})^{\frac{1}{2\la}} \, m_x=m_x\, x_Y=m_Y\,.
\eeq
By (\ref{trans1}) and (\ref{trans2}),
\beq
m_x\, dx\, dt=m_Y\, dY\, dT \com{and} dx\, dt= \xi(\cos^2\frac{v}{2})^{\frac{1}{2\la}}\,dY\, dT\,,
\eeq
and, for any $t$,
\beq\label{trans3}
m_x\, dx=m_Y\, dY\,, \com{and} dx=\xi(\cos^2\frac{v}{2})^{\frac{1}{2\la}} \, dY\,.
\eeq

Then we define $u$ as a function of the original variables $(x,t)$ by
\[u(x,t)=u\bigl(Y(x,t),T(t)\bigr)\,.\]
Note the fact that
the map $x(Y,T)$ may not be one-to-one does not cause any real
difficulty. Indeed, given $(x^*,t^*)$, we can choose an arbitrary
$Y^*$ such that $x(Y^*, t^*)=x^*$ and $T^*=t^*$, then
define $u(x^*,t^*)=u(Y^*,T^*)$. 
To prove that the values
of $u$ do not depend on the choice of $Y^*$, we proceed
as follows. Assume that there are two distinct points such that
$~x(Y_1, t^*)=x(Y_2,t^*)=x^*$, which shows that $x_Y(Y,t^*)=\xi(\cos^2\frac{v}{2})^{\frac{1}{2\la}}=0$ 
for $Y\in[Y_1, Y_2]$ because of the monotonicity of $Y$ on $x$ by \eqref{trans1}. 
This shows $\cos(\frac{v}{2})=0$ for $Y\in[Y_1, Y_2]$, where recall $\xi>0$. By (\ref{semi}), $u_Y=0$ when $Y\in[Y_1, Y_2]$. Hence, we get $u(Y_1, t^*)=u(Y_2, t^*)$.

By (\ref{semi}), (\ref{trans}), (\ref{trans2}), we could retrieve (\ref{def_lgq}). 
\subsubsection{Proof of Theorem \ref{main}}
\begin{proof}
We first consider the regularity of the solution. 
For any given time $t$, by (\ref{trans2}), (\ref{trans3}) and (\ref{semi}),
\beq\label{regu1}
\textstyle\int_{x_1}^{x_2} |u_x|^{\frac{1}{\la}}\, dx
=\int_{Y_1}^{Y_2}{ \big|}\frac{u_Y}{\xi(\cos^2\frac{v}{2})^{\frac{1}{2\la}}} {\big|}^{\frac{1}{\la}}\, \xi(\cos^2\frac{v}{2})^{\frac{1}{2\la}}dY
=\int_{Y_1}^{Y_2} |\sin\frac{v}{2}|^\frac{1}{\lambda}
\,\xi\, dY<\infty,
\eeq
on any bounded interval $[x_1,x_2]\in \RR^+$. Hence, for any time $t$, the solution $u(\cdot,t)\in W_{loc}^{1,\frac{1}{\la}}(\mathbb{R}^+)$.

Finally, we prove that the function $u$ provides a weak solution of (\ref{spe}).  By
(\ref{semi}),  (\ref{trans2}), we have 
\bes
	&&	\iint_{\mathbb{R}^+\times\mathbb{R}^+} \left\{-u_{x}\, \bigl(\phi_t+f'(u)\,\phi_x)+(\la-1)\,f''(u)(u_x)^2\phi\right\}\, dx\, dt\nonumber\\
	&=&	\iint_{(Y,T)\in\Omega} \left\{-u_Y \phi_T+(\la-1)\,f''(u)\,\xi\, \sin^2\frac{v}{2}(\cos^2 \frac{v}{2})^{\frac{1}{2\la}-1}\phi\right\}\, dY\, dT\nonumber\\
	&=&	\iint_{(Y,T)\in\Omega} \left\{u_{YT}+(\la-1)\,f''(u)\,\xi\, \sin^2\frac{v}{2}(\cos^2 \frac{v}{2})^{\frac{1}{2\la}-1}\right\}\phi\, dY\, dT\nonumber\\
	&=&	\iint_{(Y,T)\in\Omega} \Big\{(\frac{1}{2}\xi\,\sin{v}\, (\cos^2{\frac{v}{2}})^{\frac{1}{2\la}-1})_T\nonumber\\
	&&\qquad\qquad\qquad+(\la-1)\,f''(u)\,\xi\, \sin^2\frac{v}{2}(\cos^2 \frac{v}{2})^{\frac{1}{2\la}-1}\Big\}\phi\, dY\, dT\nonumber\\
	&=&	0\,,\label{proof_thm}
\ees
where $\phi(x,t)\in C^1_c(\mathbb{R}^+\times\mathbb{R}^+)$\,.
\bigskip

When $0<\la\leq\frac{1}{3}$ and $\la=\frac{1}{2}$,  
\[u_t+ f'(u)\,u_x=S\in L^\infty\,.\] 
Now we show regularities along two directions where recall another one is given in \eqref{regu1}. By the Sobolev embedding theorem,  solution $u(x,t)$ is locally H\"older continuous on both $x$ and $t$ with exponent $1-\la$. So we  complete the proof of Theorem \ref{main}.
\end{proof}
\section{Second order wave equations}
For equation (\ref{vwl}), by introducing new characteristic coordinates, we get a semi-linear system when $0<\la<1$.
Using this semi-linear system, we discuss the regularity of the solutions for (\ref{vwl}).

\subsection{A semi-linear system on new coordinates}
In this section, we only consider the smooth solution, and derive a semi-linear system
from the smooth solution of (\ref{vwl}). 

First, we introduce the new coordinate. We define
\beq\label{0RS_def}
R:=u_t+c(u)u_x,\qquad S:=u_t-c(u)u_x.
\eeq
So we have
\[
u_t=\frac{R+S}{2} \quad\text{and}\quad u_x=\frac{R-S}{2c}.
\]
To keep the tradition, we still use $R$ here to denote the gradient variables. 
By (\ref{vwl}), we have
\beq\label{0RS_eq}
\left\{\begin{array}{rcl}
R_t-cR_x&=&\frac{c'}{4c}\big[2\la R^2+(2\la-2)S^2-2(2\la-1) RS\big]\vspace{.2cm}\\
S_t+cS_x&=&\frac{c'}{4c}\big[2\la S^2+(2\la-2) R^2-2(2\la-1) RS\big]
\end{array}\right.
\eeq
where $c'=\frac{d}{du}c(u)$.
We define the forward and backward characteristics passing the point $(x,t)$ as follows
\begin{equation}
\left\{
\begin{array}{ll}
\dfrac{d}{ds}x^\pm(s;x,t)=\pm c(u(s,x^\pm(s;x,t))),\\[2mm]
x^\pm|_{s=t}=x.
\end{array}\right. \label{0char}
 \end{equation}
Along the characteristics, we define the new coordinate:
\[ X :=
\int_0^{x^-(0;x,t)}[1+R^2(y,0)]^\frac{1}{2\lambda}\,dy \com{and}
Y  :=  \int_{x^+(0;x,t)}^0[1+S^2(y,0)]^\frac{1}{2\lambda}\,dy\,.
\]
This implies
\begin{equation} 
X_t-c(u)X_x=0,\quad Y_t+c(u)Y_x=0\,. \label{0XY_con}
\end{equation}
Here, again without ambiguity, we still use $Y$ as a new coordinate for \eqref{vwl}.
For any smooth function $m,$ we obtain by using (\ref{0XY_con}) that
\begin{equation}\begin{split}
&m_t+c(u)m_x=(X_t+c(u)X_x)m_X=2c(u)X_x m_X\\
&m_t-c(u)m_x=(Y_t-c(u)Y_x)m_Y=-2c(u)Y_x m_Y.
\end{split} \label{0f_XY}
 \end{equation}
 
Next, we introduce some new variables for the semi-linear system. Without ambiguity, we still use $v$ to denote an unknown in the semi-linear system in this section. So the reader can easily compare
the semi-linear systems for \eqref{spe} and \eqref{vwl}.
We define
   \begin{equation} w  :=
2\arctan R\quad \mbox{and}\quad v  :=
2\arctan S, \label{0wz_def}
\end{equation}
so we have
\beq\label{0wz_def2}
\frac{1}{1+R^2}=\cos^2\frac{w}{2},\qquad \frac{R}{1+R^2}=\frac{1}{2}\sin w, 
\eeq
and
\beq\label{0wz_def3}
\frac{1}{1+S^2}=\cos^2\frac{v}{2}, \qquad\frac{S}{1+S^2}=\frac{1}{2}\sin v. 
\eeq
Then, we define
\begin{equation} p  :=
\frac{(1+R^2)^\frac{1}{2\la}}{X_x}\quad \mbox{and}\quad q  :=
\frac{(1+S^2)^\frac{1}{2\la}}{-Y_x}. \label{0pq_def}
\end{equation}

Finally, after having the new coordinates and new variables, we calculate the equations for $u,\ w,\ v, \ p$ and $q$. By (\ref{0RS_def}),  (\ref{0RS_eq}), (\ref{0f_XY}) and (\ref{0wz_def2})$\sim$(\ref{0pq_def}), one get the equations for $u$:
\beq\label{0u_x}
u_X=\frac{u_t+cu_x}{2cX_x}=\frac{1}{2}\sin w \,(\cos^2 \frac{w}{2})^{\frac{1}{2\la}-1} p,
\eeq
\beq\label{0u_y}
u_Y=\frac{u_t-cu_x}{-2cY_x}=\frac{1}{2}\sin{v}\,(\cos^2 \frac{v}{2})^{\frac{1}{2\la}-1}   q.
\eeq
Then, from (\ref{0RS_eq}) 
\begin{eqnarray}
w_t-cw_x&=&\frac{2}{1+R^2}(R_t-cR_x)\nonumber\\
		&=& \frac{c'}{c}\frac{1}{1+R^2}\big[\la R^2+(\la-1)S^2-(2\la-1) RS\big]\nonumber
\end{eqnarray}
so by (\ref{0RS_def}),  (\ref{0RS_eq}), (\ref{0f_XY}) and (\ref{0wz_def2})$\sim$(\ref{0pq_def}), one has
\begin{eqnarray}\label{0w_y}
w_Y
&=&\frac{w_t-cw_x}{-2cY_x}\nonumber\\
&=&\frac{c'}{2c^2}\,q\,\frac{1}{(1+S^2)^\frac{1}{2\la}}\,\frac{1}{1+R^2}\,\big[\la R^2+(\la-1)S^2-(2\la-1) RS\big]\nonumber\\
	&=&\frac{c'}{2c^2}\,q\, (\cos^2 \frac{v}{2})^{\frac{1}{2\la}-1}
	\big[ \la \sin^2\frac{w}{2} \cos^2 \frac{v}{2}+ (\la-1)\sin^2 \frac{v}{2} \cos^2 \frac{w}{2}\nonumber\\
	&&\qquad\qquad\qquad\qquad
	-\frac{2\la-1}{4}\sin w \sin v \big]
		\,.
\end{eqnarray}
Similarly, one has 
\begin{eqnarray}
v_t+cv_x&=&\frac{2}{1+S^2}(S_t+cS_x)\nonumber\\
		&=&\frac{c'}{c}\frac{1}{1+S^2}\big[\la S^2+(\la-1) R^2-(2\la-1) RS\big]\nonumber
\end{eqnarray}
and
\begin{eqnarray}\label{0v_x}
v_X
&=&\frac{v_t+cv_x}{2cX_x}\nonumber\\
&=&\frac{c'}{2c^2}\,p\,\frac{1}{(1+R^2)^\frac{1}{2\la}}\,\frac{1}{1+R^2}\,\big[\la S^2+(\la-1)R^2-(2\la-1) RS\big]\nonumber\\
	&=&\frac{c'}{2c^2}\,p\, (\cos^2 \frac{w}{2})^{\frac{1}{2\la}-1}
	\big[ \la \sin^2\frac{v}{2} \cos^2 \frac{w}{2}+ (\la-1)\sin^2 \frac{w}{2} \cos^2 \frac{v}{2}\nonumber\\
	&&\qquad\qquad\qquad\qquad
	-\frac{2\la-1}{4}\sin w \sin v \big]
		\,.
\end{eqnarray}
Lastly, we derive equations for $p$ and $q$. 
By  (\ref{0XY_con}), 
\[
X_{xt}-cX_{xx}=(X_t-cX_x)_x+c' u_x X_x=\frac{c'}{2c}(R-S)X_x  \,.
\]
So, by (\ref{0RS_def}),  (\ref{0RS_eq}), (\ref{0f_XY}) and (\ref{0wz_def2})$\sim$(\ref{0pq_def}), one has
\begin{eqnarray}
p_t-cp_x&=&(\frac{(1+R^2)^{\frac{1}{2\la}}}{X_x})_t-c (\frac{(1+R^2)^{\frac{1}{2\la}}}{X_x})_x\nonumber\\
		&=&\frac{(1+R^2)^{\frac{1}{2\la}-1}}{\la X_x}R[R_t-cR_x]-\frac{(1+R^2)^{\frac{1}{2\la}}}{X_x^2}(X_{xt}-cX_{xx})\nonumber\\
		&=&\frac{c'}{2c}\frac{(1+R^2)^{\frac{1}{2\la}-1}}{ X_x}\big[
			\frac{\la-1}{\la}RS^2 - \frac{\la-1}{\la} R^2S -R+S\big]\,,\nonumber
\end{eqnarray}
hence
\begin{eqnarray}\label{0p_x}
p_Y
&=&\frac{p_t-cp_x}{-2cY_x}\nonumber\\
&=&\frac{c'}{4c^2}\,pq\,\frac{1}{(1+S^2)^\frac{1}{2\la}}\,\frac{1}{1+R^2}\,\big[
			\frac{\la-1}{\la}RS^2 - \frac{\la-1}{\la} R^2S -R+S\big]\nonumber\\
	&=&\frac{c' }{8c^2}\,pq\,(\cos ^2\frac{v}{2})^{\frac{1}{2\la}-1}\big[ 
	\frac{\la-1}{\la}\sin w\,\sin^2 \frac{v}{2}- \frac{\la-1}{\la}\sin v\,\sin^2 \frac{w}{2}\nonumber\\
	&&\qquad -\sin w\,\cos^2 \frac{v}{2}+ \sin v\,\cos^2 \frac{w}{2}\big]\,.
\end{eqnarray}
Similarly, by  (\ref{0XY_con}), 
\[
Y_{xt}+cY_{xx}=(Y_t+cY_x)_x-c' u_x Y_x=-\frac{c'}{2c}(R-S)Y_x   \,.
\]
Then, by (\ref{0RS_def}),  (\ref{0RS_eq}), (\ref{0f_XY}) and (\ref{0wz_def2})$\sim$(\ref{0pq_def}), one has
\begin{eqnarray}
q_t+cq_x&=&(\frac{(1+S^2)^{\frac{1}{2\la}}}{-Y_x})_t+c (\frac{(1+S^2)^{\frac{1}{2\la}}}{-Y_x})_x\nonumber\\
		&=&\frac{(1+S^2)^{\frac{1}{2\la}-1}}{-\la Y_x}S[S_t+cS_x]+\frac{(1+S^2)^{\frac{1}{2\la}}}{Y_x^2}(Y_{xt}+cY_{xx})\nonumber\\
		&=&-\frac{c'}{2c}\frac{(1+S^2)^{\frac{1}{2\la}-1}}{ -Y_x}\big[
			\frac{\la-1}{\la}RS^2 - \frac{\la-1}{\la} R^2S -R+S\big]\,,\nonumber
\end{eqnarray}
hence
\begin{eqnarray}\label{0q_x}
q_X
&=&\frac{q_t+cq_x}{2cX_x}\nonumber\\
&=&-\frac{c'}{4c^2}\,pq\,\frac{1}{(1+R^2)^\frac{1}{2\la}}\,\frac{1}{1+S^2}\,\big[
			\frac{\la-1}{\la}RS^2 - \frac{\la-1}{\la} R^2S -R+S\big]\nonumber\\
	&=&-\frac{c' }{8c^2}\,pq\,(\cos^2 \frac{w}{2})^{\frac{1}{2\la}-1}\big[ 
	\frac{\la-1}{\la}\sin w\,\sin^2 \frac{v}{2}- \frac{\la-1}{\la}\sin v\,\sin^2 \frac{w}{2}\nonumber\\
	&&\qquad\qquad\qquad\qquad\qquad -\sin w\,\cos^2 \frac{v}{2}+ \sin v\,\cos^2 \frac{w}{2}\big]\,.
\end{eqnarray}

In conclusion, by (\ref{0u_x})$\sim$(\ref{0q_x}) we have a semi-linear system

\beq\label{0semi}
\left\{
\begin{array}{rcl}
u_X&=&\frac{1}{2}\sin w \,(\cos^2 \frac{w}{2})^{\frac{1}{2\la}-1} p\vspace{.1cm}\\
u_Y&=&\frac{1}{2}\sin{v}\,(\cos^2 \frac{v}{2})^{\frac{1}{2\la}-1}   q\vspace{.1cm}\\
w_Y&=&\frac{c'}{2c^2}\,q\, (\cos ^2\frac{v}{2})^{\frac{1}{2\la}-1}
	\big[ \la \sin^2\frac{w}{2} \cos^2 \frac{v}{2}+ (\la-1)\sin^2 \frac{v}{2} \cos^2 \frac{w}{2}\vspace{.1cm}\\
	&&\qquad\qquad\qquad\qquad
	-\frac{2\la-1}{4}\sin w \sin v \big]\vspace{.1cm}\\
v_X&=&\frac{c'}{2c^2}\,p\, (\cos^2 \frac{w}{2})^{\frac{1}{2\la}-1}
	\big[ \la \sin^2\frac{v}{2} \cos^2 \frac{w}{2}+ (\la-1)\sin^2 \frac{w}{2} \cos^2 \frac{v}{2}\vspace{.1cm}\\
	&&\qquad\qquad\qquad\qquad
	-\frac{2\la-1}{4}\sin w \sin v \big]\vspace{.1cm}\\
p_Y&=&\frac{c' }{8c^2}\,pq\,(\cos ^2\frac{v}{2})^{\frac{1}{2\la}-1}\big[ 
	\frac{\la-1}{\la}\sin w\,\sin^2 \frac{v}{2}- \frac{\la-1}{\la}\sin v\,\sin^2 \frac{w}{2}\vspace{.1cm}\\
	&&\qquad \qquad\qquad\qquad\quad-\sin w\,\cos^2 \frac{v}{2}+ \sin v\,\cos^2 \frac{w}{2}\big]\vspace{.1cm}\\
q_X	&=&-\frac{c' }{8c^2}\,pq\,(\cos^2 \frac{w}{2})^{\frac{1}{2\la}-1}\big[ 
	\frac{\la-1}{\la}\sin w\,\sin^2 \frac{v}{2}- \frac{\la-1}{\la}\sin v\,\sin^2 \frac{w}{2}\vspace{.1cm}\\
	&&\qquad\qquad \qquad\qquad\quad-\sin w\,\cos^2 \frac{v}{2}+ \sin v\,\cos^2 \frac{w}{2}\big]\,.
\end{array}
\right.
\eeq
%
\subsection{Regularity of solutions when $\la\in(0,\frac{1}{3}]$}
In the rest of this paper, we restrict our consideration on the case when $\la\in(0,\frac{1}{3}]$.
Recall that $c(u)$ is assumed to be uniformly positive and bounded.

If we first suppose that there is a $L^\infty$ bound on $p$ and $q$ on any bounded set of the $(X,Y)$-plane (on $t\geq0$ part), then the right hand side of semi-linear system (\ref{0semi}) is Lipschitz  on $u,\ v,\ w,\ p$ and $q$ when $\la\in(0,\frac{1}{3}]$ on any bounded set of the $(X,Y)$-plane (on $t\geq0$ part). 

Hence,  using similar method as those used for \eqref{spe} in Theorem \ref{main} and for \eqref{vw} in \cite{BZ}, we could prove the global existence of solution for semi-linear system (\ref{0semi}) then recover a weak solution for 
\eqref{vwl}, where $u_t\pm c(u)u_x$ are in $\textstyle W^{1,\frac{1}{\la}}$, because we assume that $p$ and $q$ are bounded.
Hence, the solution $u$ is H\"older continuous with exponent $1-\la$. We leave the details to the reader
and also refer the reader to \cite{BZ}.

Unfortunately, till now, we still have no method to bound  $p$ and $q$, hence the existence and regularity of the solution is still a conjecture and only can be supported by numerical experiments which will be shown later. For \eqref{vw}, the bounds on $p$ and $q$ are found in \cite{BZ} by exploring energy law. But equation \eqref{vwl} with $\la\neq \frac{1}{2}$ is not endowed with such an energy law.
Instead we have an equation
\beq
(\cos^2 \frac{w}{2})^{\frac{1}{2\la}-1} p_Y+(\cos^2\frac{v}{2})^{\frac{1}{2\la}-1} q_X=0\,.
\eeq

On the other hand, system \eqref{0semi} still can help us analyze the original equation \eqref{vwl}. 
As we discussed, to show whether $u$ is H\"older continuous with exponent $1-\la$,
by numerical  experiments, we only have to test whether $p$ and $q$ remain uniformly positive and bounded
when $v$ or $w$ attains $\pi$, i.e. when gradient blowup happens. 

Recall that the region $(x,t)\in{\mathbb R}\times {\mathbb R}^+$ is transformed to the region $\{Y\geq\phi(X)\}$ in $(X,Y)$ plane, where $t=0$ is corresponding to the curve $Y=\phi(X)$ with $\phi(X)$ strictly decreasing on $X$. Then 
 we do some numerical experiments on the bounded domain of 
 \[
\bar\Omega_r=\left\{ (X,Y);\ Y\geq \phi(X),\quad X\leq r,\quad Y\leq r  \right\}.
\] 
We can iteratively find the solution by initially set
\beq\label{initial}
\left\{
\begin{array}{rcl}
u_0(X,Y) &=& u(\phi^{-1}(Y), Y) \vspace{.1cm} \\
w_0(X,Y)&=&w(X, \phi(X))\vspace{.1cm}\\
v_0(X,Y)&=&v(\phi^{-1}(Y), Y)\vspace{.1cm}\\
p_0(X,Y)&=&p(X,\phi(X))\vspace{.1cm}\\
q_0(X,Y)&=&q(\phi^{-1}(Y),Y)\,.
\end{array}
\right.
\eeq
with $u(\phi^{-1}(Y), Y) $ given by $u(x,t=0) = sech(x)$. We also set $u_t(x,t=0) = u_x(x,t=0)$ initially and use $c(u) = \sqrt{\cos^2(u)+1}$. Then we do the iteration according to the following integral equations coming from
\eqref{0semi}.
\beq\label{iterate}
\small{\left\{
\begin{array}{l}
\textstyle u_{n+1}(X,Y) = u_n(\phi^{-1}(Y), Y)+\int_{\phi^{-1}(Y)}^X \frac{1}{2}\sin w_{n} \,(\cos^2 \frac{w_{n}}{2})^{\frac{1}{2\la}-1} p_{n} dX \vspace{.1cm} \\
\textstyle w_{n+1}(X,Y)=w_n(X, \phi(X))+\int_{\phi(X)}^Y \frac{c'}{2c^2}\,q\, (\cos ^2\frac{v_n}{2})^{\frac{1}{2\la}-1}
	\big[ \la \sin^2\frac{w_n}{2} \cos^2 \frac{v_n}{2}\vspace{.1cm} \\
	\qquad\qquad\qquad
	 \textstyle+(\la-1)\sin^2 \frac{v_n}{2} \cos^2 \frac{w_n}{2}-\frac{2\la-1}{4}\sin w \sin v \big]dY\vspace{.1cm}\\
\textstyle v_{n+1}(X,Y)=v_n(\phi^{-1}(Y), Y)+\int_{\phi^{-1}(Y)}^X \frac{c'}{2c^2}\,p_n\, (\cos^2 \frac{w_n}{2})^{\frac{1}{2\la}-1}
	\big[ \la \sin^2\frac{v_n}{2} \cos^2 \frac{w_n}{2}\vspace{.1cm}\\
	\qquad\qquad\qquad
	\textstyle+ (\la-1)\sin^2 \frac{w_n}{2} \cos^2 \frac{v_n}{2}-\frac{2\la-1}{4}\sin w_n \sin v_n \big] dX\vspace{.1cm}\\
\textstyle p_{n+1}(X,Y)=p_n(X,\phi(X))+ \int_{\phi(X)}^Y \frac{c' }{8c^2}\,p_n q_n\,(\cos ^2\frac{v_n}{2})^{\frac{1}{2\la}-1}\big[ 
	\frac{\la-1}{\la}\sin w_n\,\sin^2 \frac{v_n}{2}\vspace{.1cm}\\
	\textstyle\qquad\qquad\qquad - \frac{\la-1}{\la}\sin v_n\,\sin^2 \frac{w_n}{2} -\sin w_n\,\cos^2 \frac{v_n}{2}+ \sin v_n\,\cos^2 \frac{w_n}{2}\big] dY\vspace{.1cm}\\
\textstyle q_{n+1}(X,Y)=q_n(\phi^{-1}(Y),Y)-\int_{\phi^{-1}(Y)}^X \frac{c' }{8c^2}\,p_n q_n\,(\cos^2 \frac{w_n}{2})^{\frac{1}{2\la}-1}\big[ 
	\frac{\la-1}{\la}\sin w_n\sin^2 \frac{v_n}{2}\vspace{.1cm}\\
	\textstyle\qquad\qquad\qquad - \frac{\la-1}{\la}\sin v_n\,\sin^2 \frac{w_n}{2}-\sin w_n\,\cos^2 \frac{v_n}{2}+ \sin v_n\,\cos^2 \frac{w_n}{2}\big] dX\,.
\end{array}
\right.}
\eeq

\begin{figure}[htb]
\centering
\includegraphics[width=0.45\textwidth]{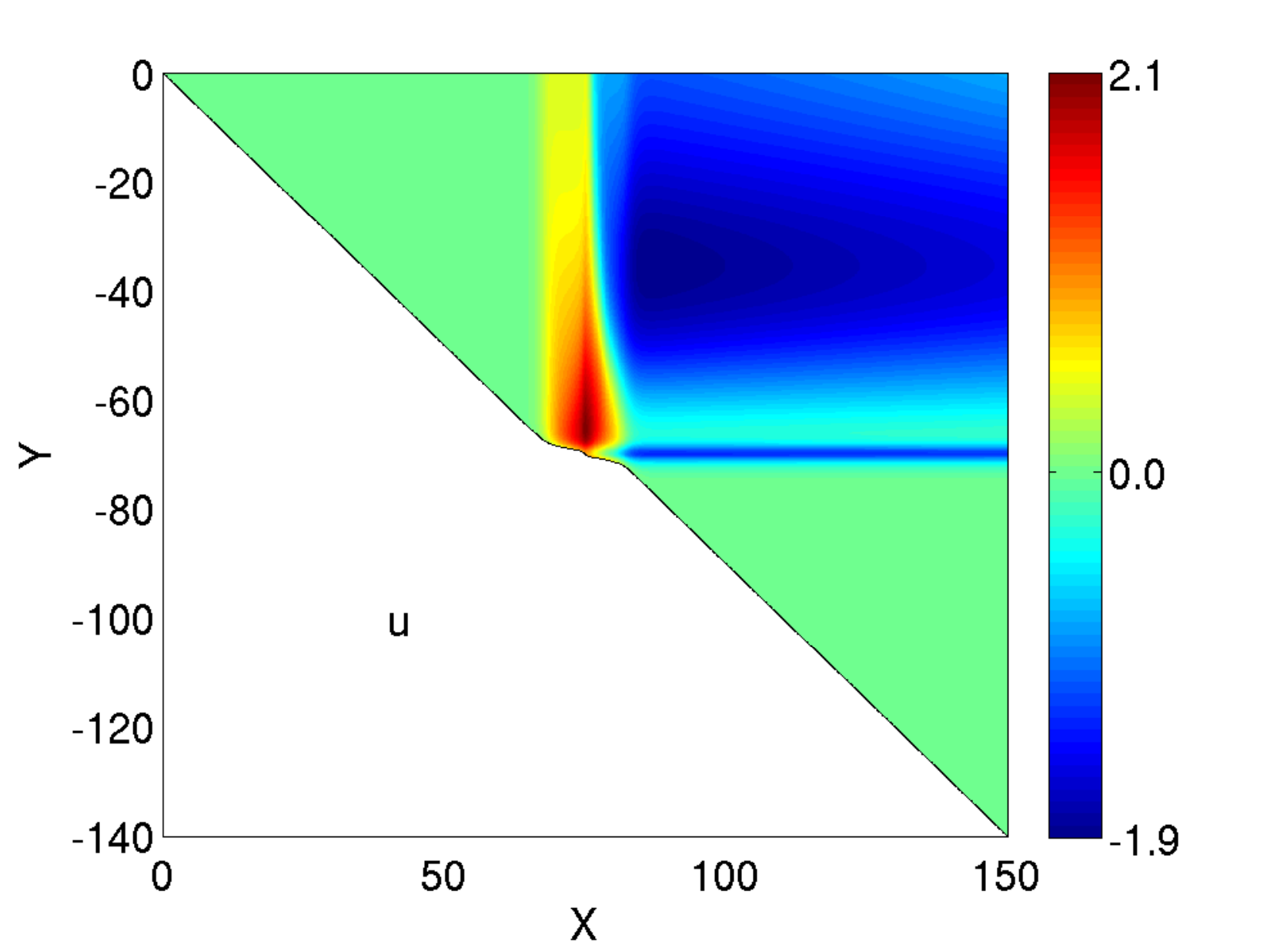}
\includegraphics[width=0.45\textwidth]{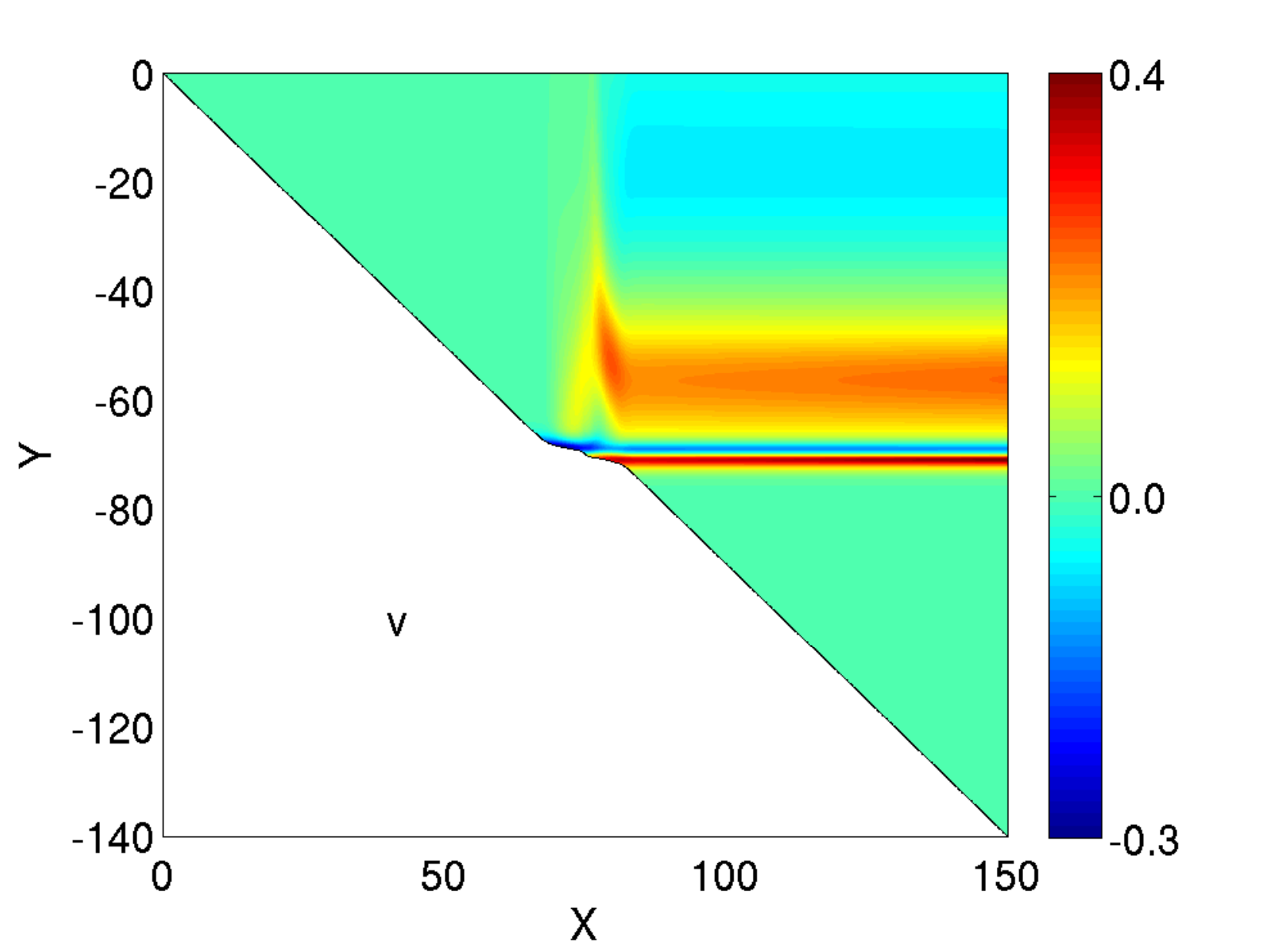}
\includegraphics[width=0.45\textwidth]{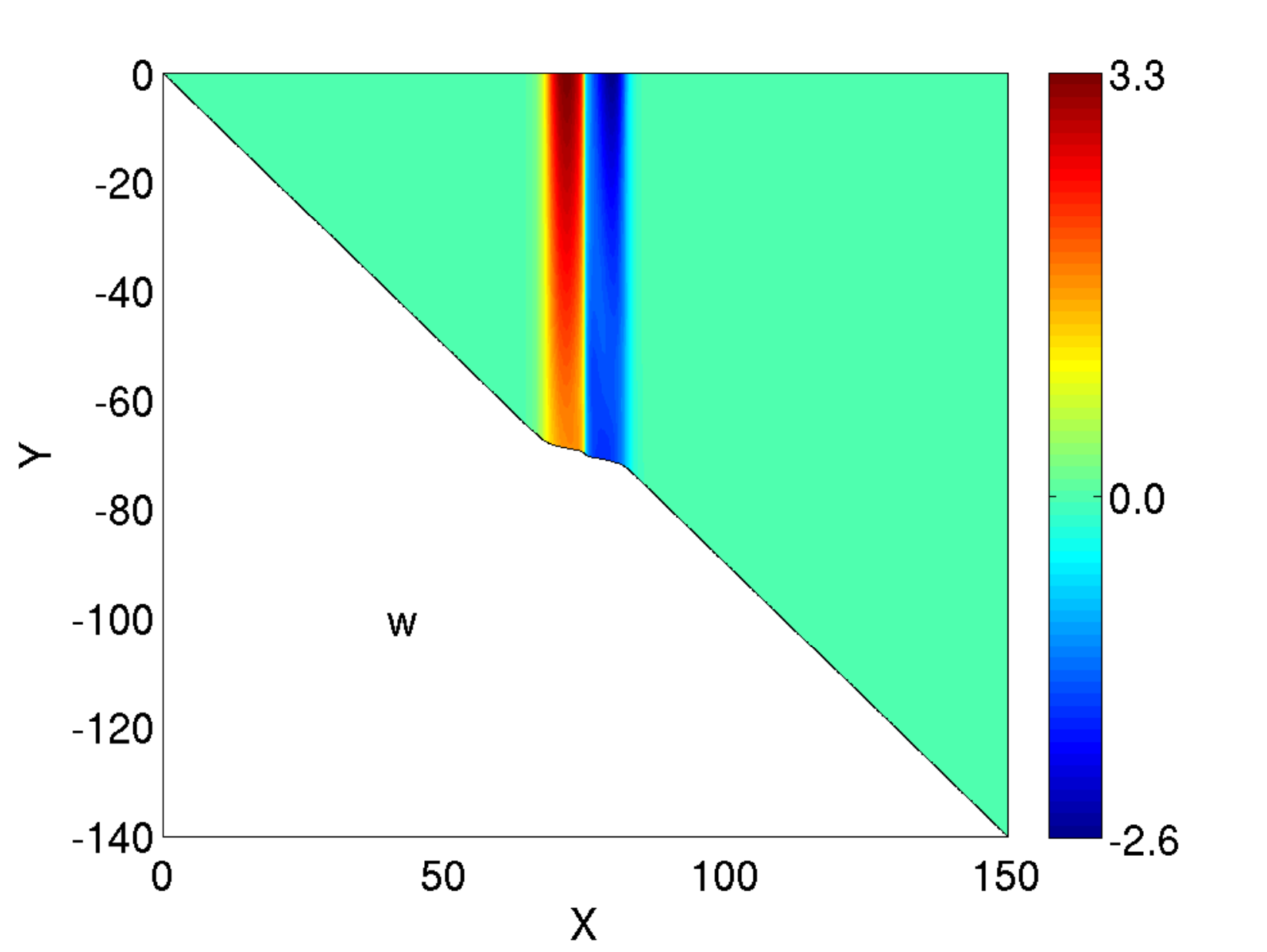}
\includegraphics[width=0.45\textwidth]{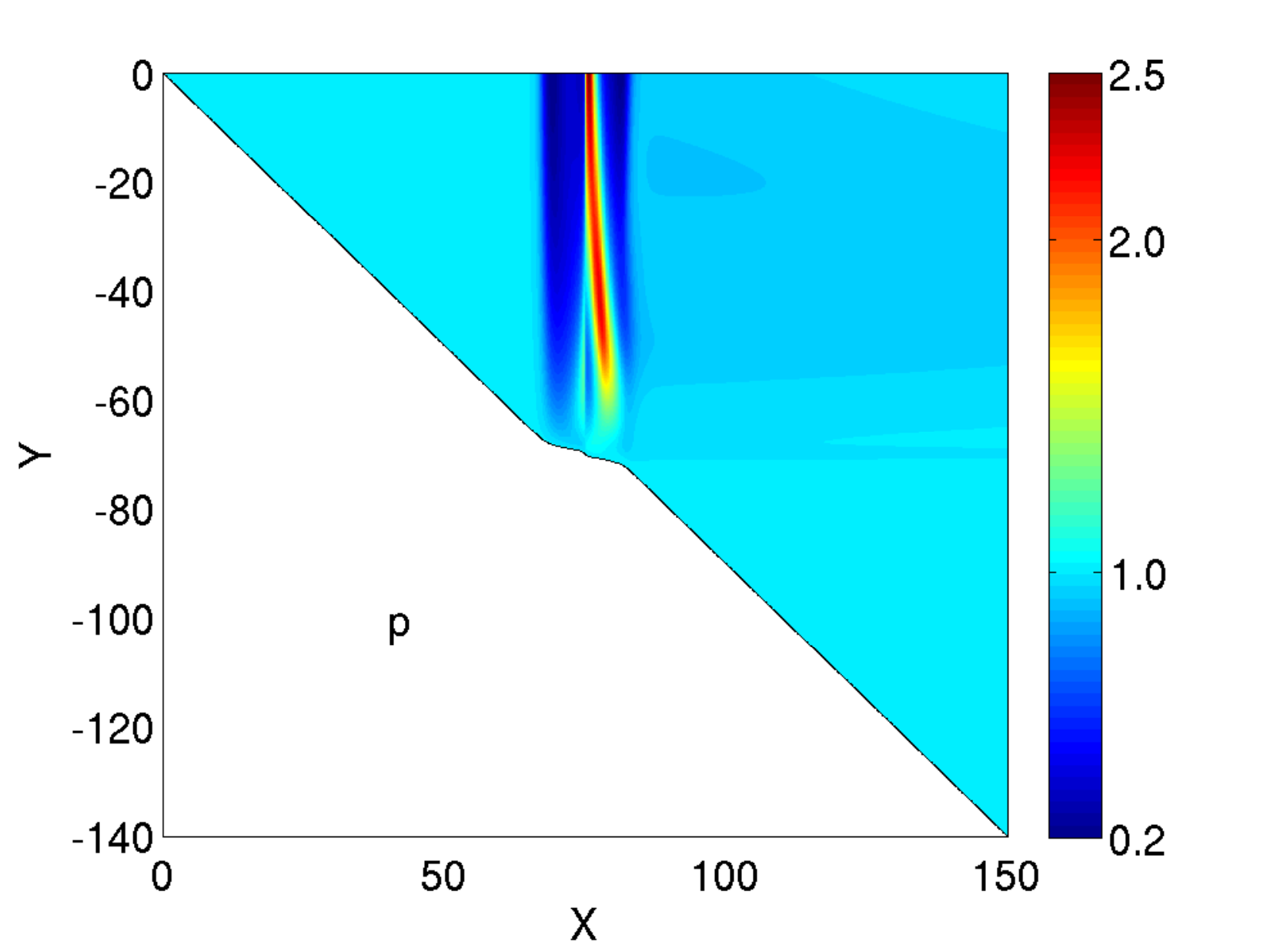}
\includegraphics[width=0.45\textwidth]{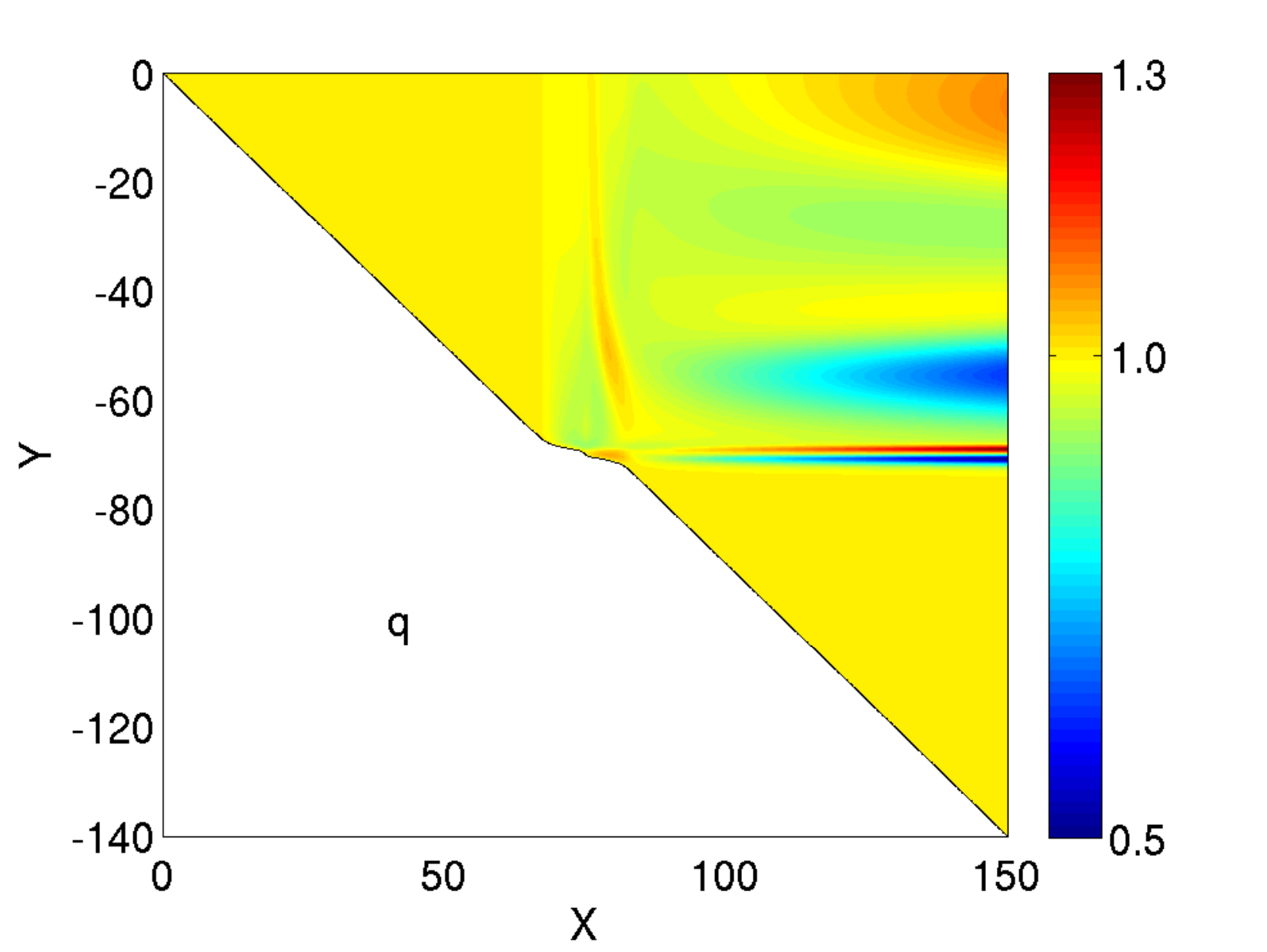}
\caption{$u,\ v,\ w,\ p,\ q$ in $(X,Y)$ coordinate for $\lambda=\frac{1}{4}$. In this example, $w$ attains $\pi$ in finite time, but $p$ and $q$ are both uniformly bounded and positive.}
\label{l_1o4}
\end{figure}

In a numeric experiment shown in Figure {\ref{l_1o4}} when $\la=1/4$, we see that $p$ and $q$ remain uniformly positive and bounded when $w$ attains $\pi$, which strongly indicates that we could construct a solution $u(x,t)$ which is H\"older continuous with exponent $1-\la$, including finite time blowup ($w$ attains $\pi$). 

In fact, on the bounded set $\bar\Omega_r$ if $p$ and $q$ have positive upper and lower bounds, 
then equation (\ref{iterate}) gives a contract mapping in 
\[
	\Lambda_r := \left\{ f : \bar\Omega_r \mapsto \mathbb{R};\  \| f \|_*:= 
	{\text{ess}}\sup_{(X,Y)\in \bar\Omega_r} e^{-\kappa(|X|+|Y|)} |f(X,Y)|<\infty \right\},
\]
with sufficiently large $\kappa$. Hence $u_n,\ v_n,\ w_n,\ p_n,\ q_n$ converge to a solution of equation (\ref{0semi}). Then we could construct a solution $u(x,t)$ for \eqref{vwl} which is H\"older continuous with exponent $1-\la$,
through similar inverse transformation as Theorem \ref{main} and the one in \cite{BZ} for \eqref{vw}.

\begin{figure}[htb]
\centering
\includegraphics[width=0.45\textwidth]{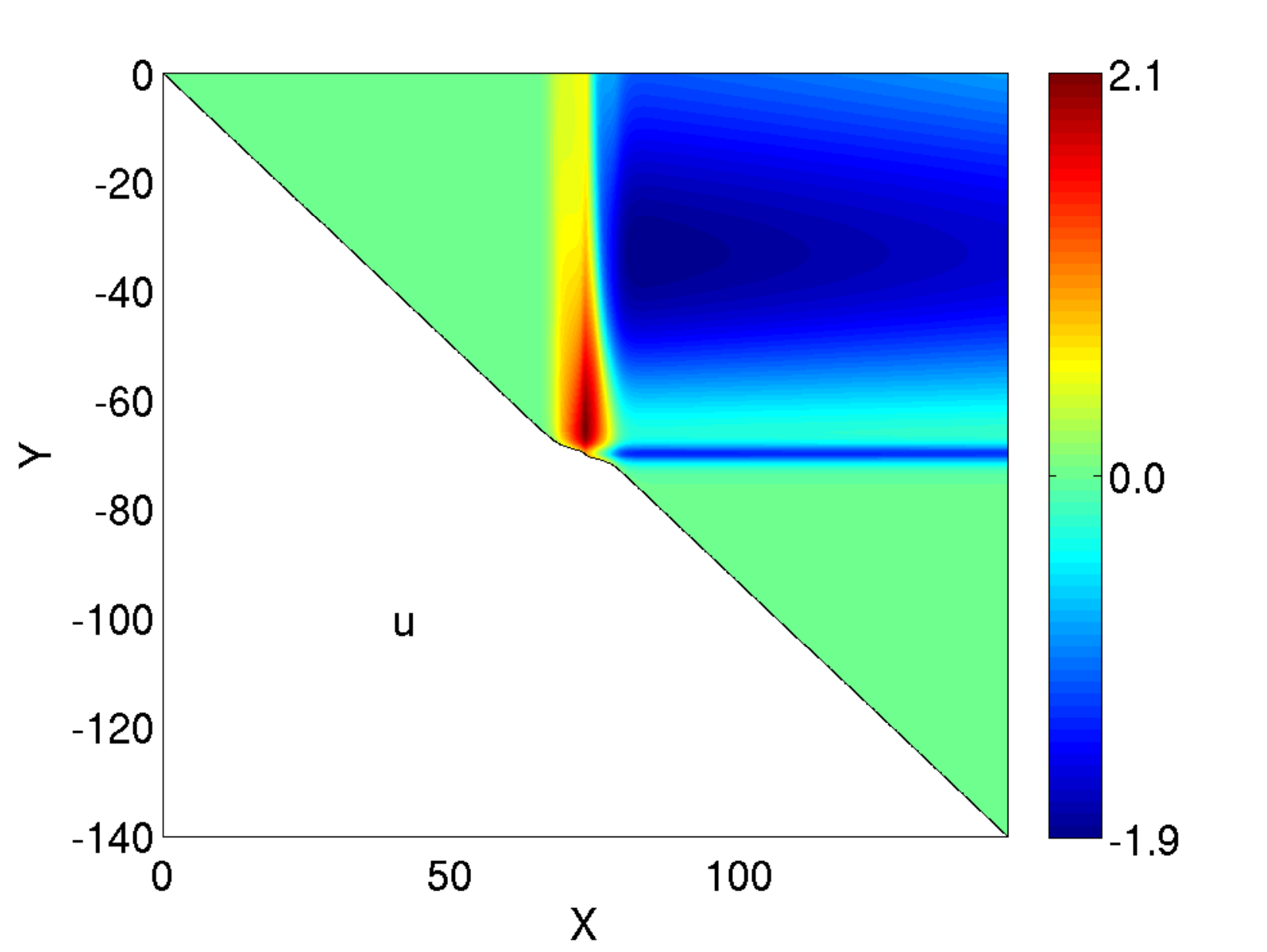}
\includegraphics[width=0.45\textwidth]{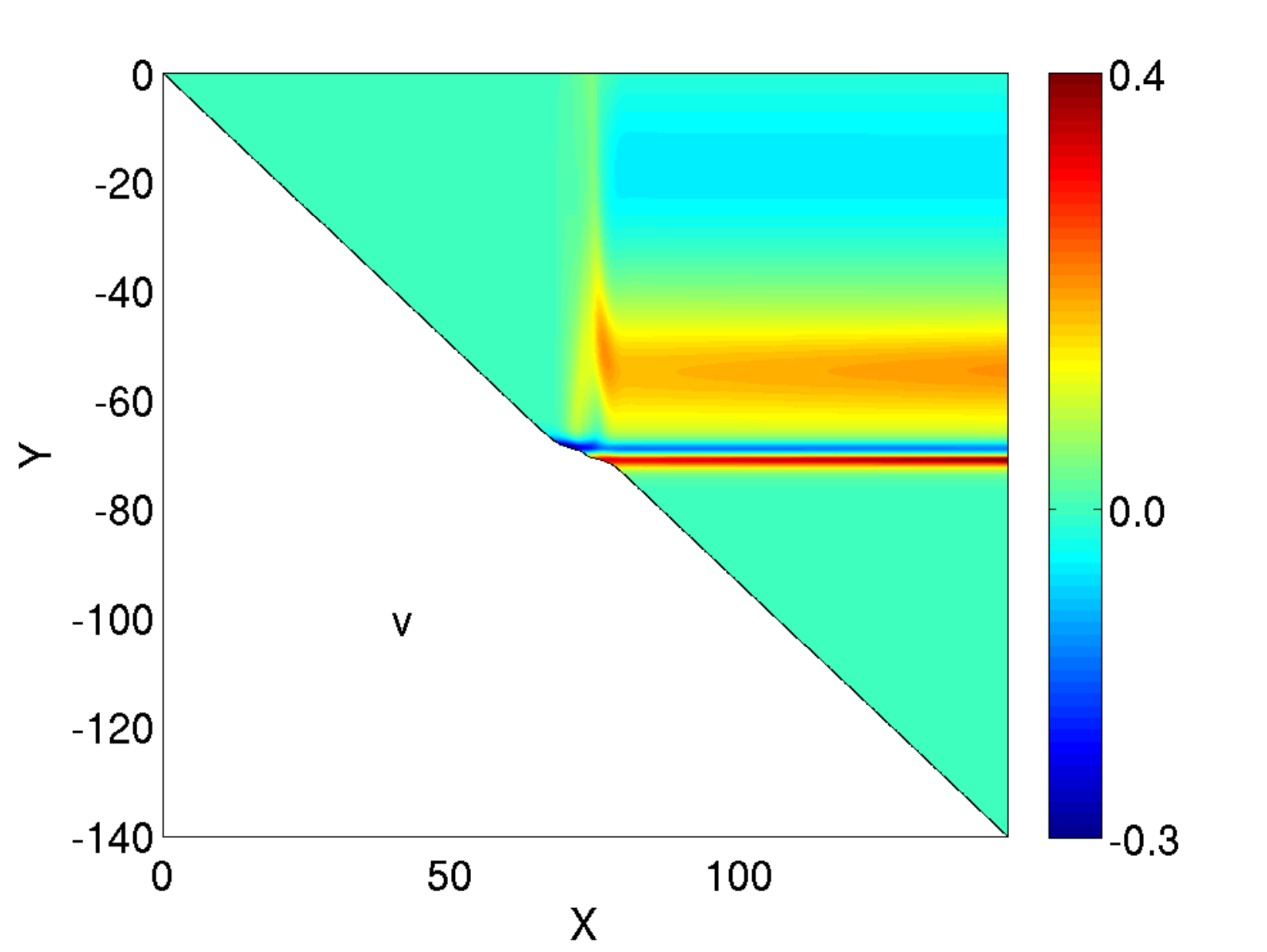}
\includegraphics[width=0.45\textwidth]{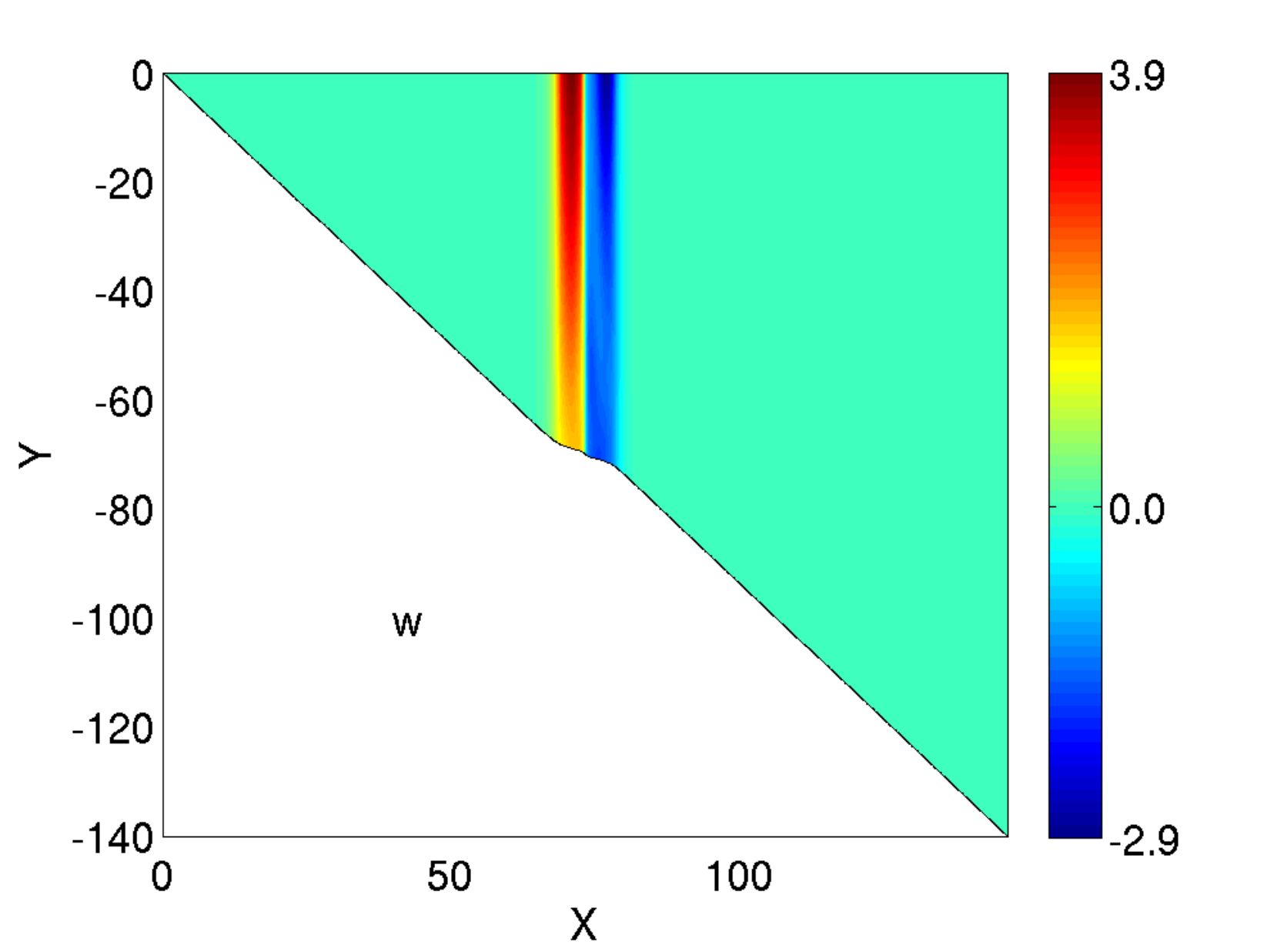}
\includegraphics[width=0.45\textwidth]{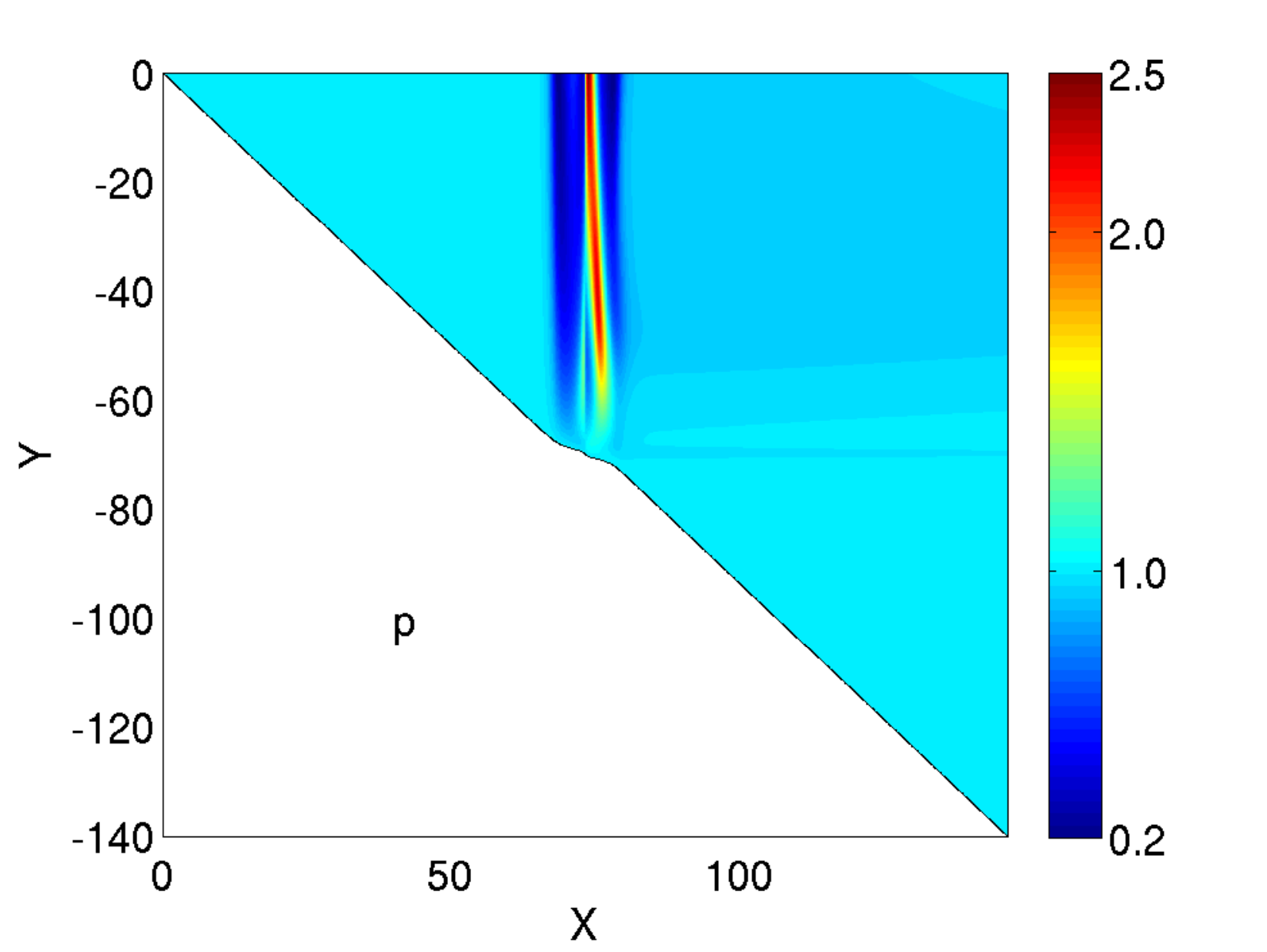}
\includegraphics[width=0.45\textwidth]{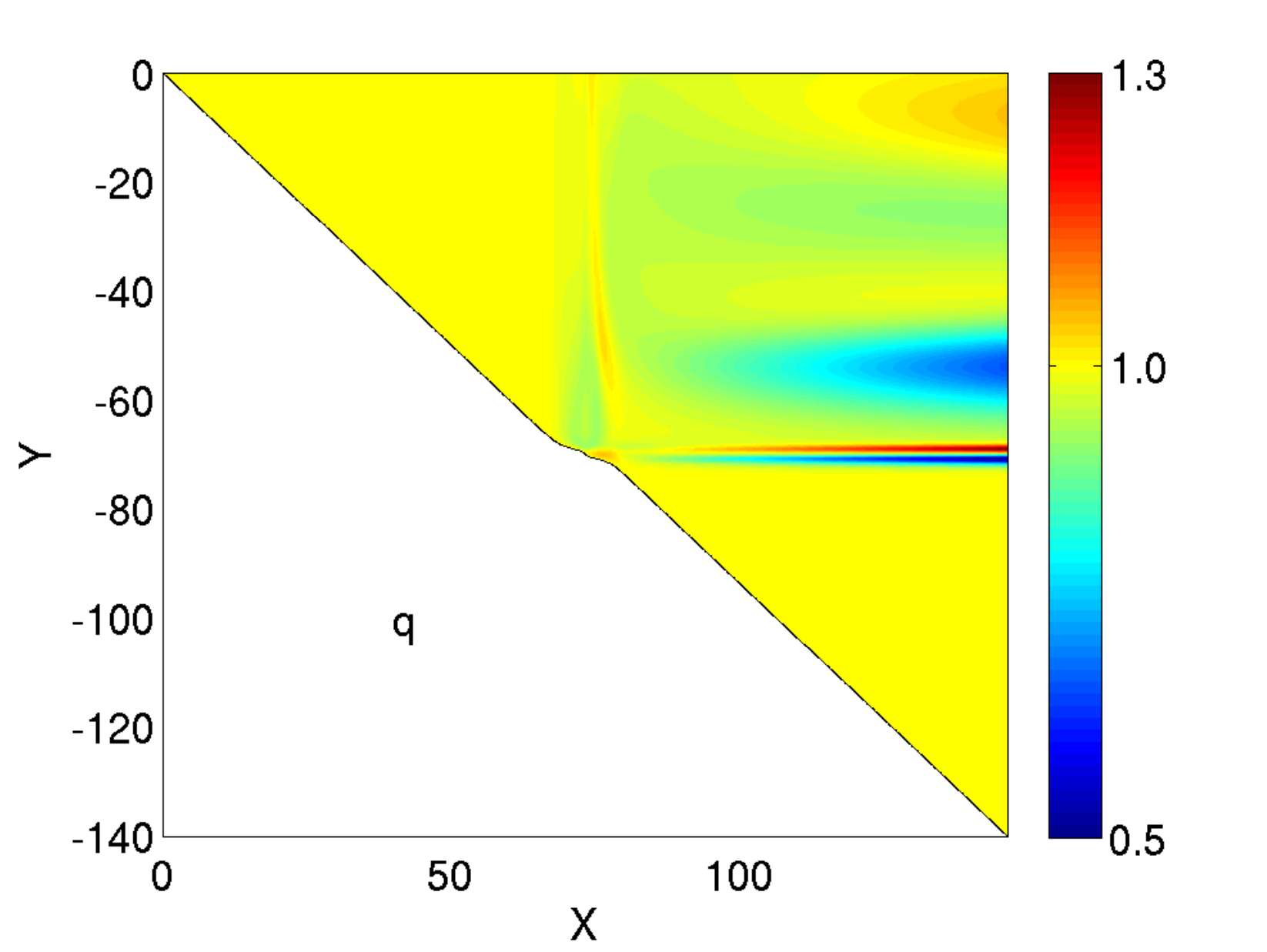}
\caption{$u,\ v,\ w,\ p,\ q$ in $(X,Y)$ coordinate for $\lambda=\frac{1}{3}$. In this example, $w$ attains $\pi$ in finite time, but $p$ and $q$ are both uniformly bounded and positive.}
\label{l_1o3}
\end{figure}

Similar experiment has been done for $\la=1/3$, which is shown in  Figure {\ref{l_1o3}}. In this case, we can draw a same conclusion as case $\la=1/4$.

The advantage of doing numeric analysis on the semi-linear system \eqref{0semi} instead of analyzing \eqref{vwl} directly is that we change a problem for H\"older continuity for $u$ when gradient blowup happens into a problem for $L^\infty$ bounds of $p$ and $q$ when $v$ or $w$ attains $\pi$.  Clearly, it is much easier to check the latter one than the first one by numerical methods. The semi-linear system is also easier to cope with than the quasi-linear system.

\begin{thebibliography}{88}
\bibitem{BC} A.~Bressan and A.~Constantin, \emph{Global solutions of the Hunter-Saxton equation},
SIAM J. Math. Anal., {\bf 37} (2005), 996--1026.

\bibitem{BC2} A.~Bressan and A.~Constantin, \emph{Global conservative solutions to the Camassa-Holm equation}, Arch. Rat. Mech. Anal. {\bf183} (2007), 215-239.

\bibitem{BHR} A.~Bressan, H.~Holden and X.~Raynaud, 
   \emph{Lipschitz metric for the Hunter-Saxton equation}, J. Math. Pures Appl. (9) {\bf 94}:1 (2010), 68Ð92. 

\bibitem{BZ}A.~Bressan and Y.~Zheng,  \emph{Conservative solutions to a
nonlinear variational wave equation},
Comm. Math. Phys., {\bf 266} (2006), 471--497.

\bibitem{BZZ} Alberto Bressan, Ping Zhang, and Yuxi Zheng: \emph{Asymptotic
variational wave equations},  Arch. Ration. Mech. Anal.,
\textbf{183} (2007),  163--185.

\bibitem{G3}  G.~Chen,
{\em Formation of singularity and smooth wave propagation for the
non-isentropic compressible Euler
equations,} { J. Hyperbolic Differ. Equ.}, \textbf{8}:4 (2011), 671--690.

\bibitem{CPZ} G.~Chen, R.~Pan and S.~Zhu,
{\em Singularity formation for compressible Euler equations,} { submitted}, available at arXiv:1408.6775.


\bibitem{G5} G.~Chen and R.~Young, {\em Smooth waves and gradient blowup for the inhomogeneous 
wave equations,} {J. Differential Equations}, {\bf 252}:3 (2012), 2580--2595.

\bibitem{G8} G.~Chen, R.~Young and Q.~Zhang, {\em Shock formation in the compressible 
Euler equations and related systems,} { J. Hyperbolic Differ. Equ.}, \textbf{10}:1 (2013), 149--172.

\bibitem{CZZ12} G.~Chen, P.~Zhang and Y.~Zheng, Energy Conservative Solutions to a Nonlinear Wave
System of Nematic Liquid Crystals,  {\it Comm. Pure Appl. Anal.}, \textbf{12}:3 (2013), 1445--1468.

\bibitem{GCZ} G.~Chen and Y.~Zheng, {\em Existence and singularity to a wave system of 
nematic liquid crystals,}  J. Math. Anal. Appl., \textbf{398} (2013), 170--188.

\bibitem{Dafermos} C.~M.~Dafermos, {Hyperbolic Conservations laws in
Continuum Physics (third edition)}, \emph{Springer-Verlag,} Heidelberg, 2010.

\bibitem{HR} H. Holden and X. Raynaud, 
 {\em Global semigroup of conservative solutions of the nonlinear variational wave equation. }  {Arch. Ration. Mech. Anal.} {\bf 201}:3 (2011), 871Ð964.

\bibitem{ghz}  R.~Glassey, J.~Hunter and Y.~Zheng,
Singularities of a variational wave equation, {\it J. Differential Equations}, 
{\bf 129} (1996), 49--78.

\bibitem{ghz2}   R.~Glassey, J.~Hunter and Y.~Zheng,
Singularities and Oscillations in a Nonlinear Variational Wave Equation, {\it Singularities and Oscillations, The IMA Volumes in Mathematics and its Applications,} {\bf 91} (1997), 37-60.

\bibitem{HS}
J.~K.~Hunter and R.~H.~Saxton, \emph{Dynamics of director fields}, 
SIAM J. Appl. Math., {\bf 51} (1991), 1498-1521.

\bibitem{HZ95a}  J.~K.~Hunter and Yuxi Zheng,
\emph{On a Nonlinear Hyperbolic Variational Equation: I. Global Existence
of Weak Solutions},  Arch. Rat. Mech. Anal., {\bf 129} (1995),
305--353.

\bibitem{HZ95b}  J.~K.~Hunter and Yuxi Zheng,
\emph{On a Nonlinear Hyperbolic Variational Equation: II. The Zero Viscosity and Dispersion Limits}, 
 Arch. Rat. Mech. Anal., {\bf 129} (1995),
355--383.


\bibitem{lax0} P.~Lax, \emph{ Development of singularities of solutions
of nonlinear hyperbolic partial differential equations,} J.
Math. Physics, {\bf 5}:5 (1964) 611-614.

\bibitem{LG}P.~Lax and J.~Glimm,
Decay of solutions of systems of nonlinear hyperbolic conservation laws,
\emph{Memoirs of the American Mathematical Society,} {\bf 101}, American Mathematical Society, Providence, R.I. 1970. 


\bibitem{lindblad}H. Lindblad,  \emph{Global solutions of nonlinear wave equations.} Comm. Pure Appl. Math., {\bf 45}:9 (1992), 1063Ð1096.

\bibitem{ZZ03} 
\newblock Ping Zhang and Yuxi Zheng,  
\newblock \emph{Weak solutions to a nonlinear variational wave equation}, 
\newblock Arch. Ration. Mech. Anal., \textbf{166} (2003), 303--319.

\bibitem{ZZ05a} 
\newblock Ping Zhang and Yuxi Zheng,
\newblock \emph{Weak solutions to a nonlinear variational wave equation with general data}, 
\newblock  Ann. I. H. Poincar\'e, \textbf{22} (2005), 207--226.

\bibitem{ZZ10} Ping Zhang and Yuxi Zheng, \emph{Conservative
solutions to a  system of variational wave equations of nematic
liquid crystals},  Arch. Ration. Mech. Anal., {\bf 195} (2010), 701-727.

\bibitem{ZZ11} Ping Zhang and Yuxi Zheng, \emph{Energy Conservative Solutions to a One-Dimensional Full Variational Wave
System},  Comm. Pure Appl. Math., {\bf 55} (2012), 582-632.

\end{thebibliography}
\end{document}